\documentclass[12pt,openbib]{article}%
\usepackage{amsmath}
\usepackage{amsfonts}
\usepackage{amssymb}
\usepackage{graphicx}
\usepackage{natbib}%
\usepackage{setspace}
\newtheorem{theorem}{Theorem}

\newenvironment{proof}[1][Proof]{\noindent\textbf{#1.} }{\ \rule{0.5em}{0.5em}}

\newcommand{\jmdouble}{\relax}
\jmdouble
\begin{document}

\title{Fast computation by block permanents of cumulative distribution functions of
order statistics from several populations\footnotemark[1]}
\author{D. H. Glueck\footnotemark[2]{\ }\footnotemark[4]
\and A. Karimpour-Fard\footnotemark[2]
\and J. Mandel\footnotemark[2]
\and L. Hunter\footnotemark[2]
\and K. E. Muller\footnotemark[3]}
\date{\today}
\maketitle

\footnotetext[1]{Deborah H. Glueck is Assistant Professor, Department of
Preventive Medicine and Biometrics, University of Colorado at Denver and
Health Sciences Center, Campus Box B119, 4200 East Ninth Avenue, Denver,
Colorado 80262 (e-mail: Deborah.Glueck@uchsc.edu). Anis Karimpour-Fard is a
graduate student in Bioinformatics, Department of Preventive Medicine and
Biometrics, University of Colorado at Denver and Health Sciences Center,
Campus Box B119, 4200 East Ninth Avenue, Denver, Colorado 80262 (e-mail: Anis
Karimpour-Fard@uchsc.edu). Jan Mandel is Professor, Department of Mathematics,
Adjunct Professor, Department of Computer Science, and Director of the Center
for Computational Mathematics, University of Colorado at Denver and Health
Sciences Center, Campus Box 170, Denver, Colorado 80217-3364
(e-mail:Jan.Mandel@cudenver.edu). Larry Hunter is Associate Professor of
Biology, Computer Science, Pharmacology, and Preventive Medicine and
Biometrics, and Director of the Center for Computational Pharmacology,
University of Colorado at Denver and Health Sciences Center, PO Box 6511, MS
8303, Aurora, CO 80045-0511 (e-mail: Larry.Hunter@uchsc.edu). Keith E. Muller
is Professor and Director of the Division of Biostatistics, Department of
Epidemiology and Health Policy Research, University of Florida, 1329 SW 16th
Street Room 5125, PO Box 100177 Gainesville, FL 32610-0177
(e-mail:Keith.Muller@biostat.ufl.edu)

Glueck was supported by NCI K07CA88811.
Mandel was supported by NSF-CMS 0325314. Muller was supported by NCI P01 CA47
982-04, NCI R01 CA095749-01A1 and NIAID 9P30 AI 50410. Hunter was supported by
NIAAA 1U01 AA13524-02 and NCI 5 P30 CA46934-15.

The authors thank Professor Gary Grunwald for his helpful comments.}

\footnotetext[2]{University of Colorado at Denver and Health Sciences Center}
\footnotetext[3]{University of Florida}

\newpage

\begin{abstract}
The joint cumulative distribution function for order statistics arising from
several different populations is given in terms of the distribution function
of the populations. The computational cost of the formula in the case of two
populations is still exponential in the worst case, but it is a dramatic
improvement compared to the general formula by Bapat and Beg. In the case when
only the joint distribution function of a subset of the order statistics of
fixed size is needed, the complexity is polynomial, for the case of two populations.
\end{abstract}


\newpage

\section{INTRODUCTION}

The \citet{Benjamini-1995-CFD} procedure represents one of what has become a
rather large class of techniques in which we would like to be able to
calculate order statistics arising from several populations. The complexity of
the calculations implied by such approaches has remained a barrier to accurate
probability statements. We provide tools which greatly extend the range of
computable cases.

Order statistics obtained by sampling from two different populations occur,
e.g., when $p$-values arise from null or alternative hypotheses, from men or
women, or from two different types of cancer.

The distribution of order statistics for independent, identically distributed
random variables is well known, and appears in every basic statistics book;
for example, \citet[Chapter 4, Section
6]{Hogg-1978-IMS}. \citet{David-2003-OS} and \citet{Balakrishnan-1998-OSI}
provide a thorough review of order statistics. For identically distributed
random variables, the cumulative distribution function is concise and fast to compute.

For independent, but not identically distributed random variables, a formula
for computing the joint cumulative distribution function of the order
statistics was given by \citet{Bapat-1989-OSN}. However, this formula is
computationally intractable, because it involves an exponential number of
permanents of the size of the number of random variables. In addition, the
complexity of the computation of the permanent by the best algorithms grows
exponentially \citep[p.~499]{Knuth-1998-TAC}. Approximate algorithms for
computing the permanent
\citep{Valiant-1979-CCP,Forbert-2003-CPS,Jerrum-2004-PAA} with lower
asymptotic complexity are still not practical.

We show that the computational cost of the formula in the case of two
populations is still exponential, but is a dramatic improvement compared to
the general formula by Bapat and Beg. In the case when only the joint
distribution function of a subset of the order statistics of fixed size is
needed, we show that the complexity is polynomial, in the case of two populations.

\section{NOTATION AND PRELIMINARIES}

For an $m\times m$ matrix $\boldsymbol{A}$, with entries $a_{ij}$, the
permanent is given by \citet[p. 30]{Aitken-1939-DM}
\begin{equation}
\operatorname{per}\left[  \boldsymbol{A}\right]  =\sum_{\pi}\prod_{i=1}%
^{m}a_{i,\pi\left(  i\right)  }\text{ .} \label{eq:permanent}%
\end{equation}
where $\pi$ ranges over all permutations of $\left\{  1,2,\ldots,m\right\}  $.
Hence, the permanent is defined much like the determinant, but with all signs
positive. The permanent can be expanded by row or columns exactly like the
determinant. The computational cost of evaluating the permanent by expansion
is $O(m!)$ operations. The computational cost using the best algorithms is
exponential \citet[p. 499]{Knuth-1998-TAC}.

The following notation will be used in all theorems and proofs in this paper
without further explicit reference. $X_{i}$, $i=1,\ldots,m$ are independent
real valued random variables with cumulative distribution functions
$F_{i}\left(  x\right)  $. \ The order statistics $Y_{1},Y_{2},\ldots,Y_{m}$
are random variables defined by sorting the values of $X_{i}$. In particular,
$Y_{1}\leq Y_{2}\leq\ldots\leq Y_{m}$. The arguments of the joint cumulative
distribution function of order statistics are customarily written omitting
redundant arguments; thus let $n\,,$ $1\leq n_{1}<n_{2}<\cdots<n_{_{k}}\leq
m$, denote the indices of the remaining arguments and $y_{1}\leq y_{2}%
\leq\cdots\leq y_{k}$ their values. Finally, define the index vector
$\mathbf{i}=\left(  i_{0},i_{1},\ldots i_{k+1}\right)  $ and the summation
index set%
\begin{equation}
\mathcal{I}=\left\{  \mathbf{i:}%
\begin{array}
[c]{c}%
0=i_{0}\leq i_{1}\leq\cdots\leq i_{k}\leq i_{k+1}=m\text{, }\\
\text{and }i_{j}\geq n_{j}\text{ for all }1\leq j\leq k
\end{array}
\right\}  . \label{eq:def-I}%
\end{equation}
Writing summation over the set $\mathcal{I}$ in terms of loops is
straightforward. Using the set $\mathcal{I}$ instead of the loop in this paper
allows an insight into the structure of the method and its complexity, and it
does not tie the mathematical formulation to any particular implementation.

The joint cumulative distribution function of the set $\left\{  Y_{n_{1}%
},Y_{n_{2}},\ldots,Y_{n_{k}}\right\}  $, which is a subset of the complete set
of order statistics, is defined as%
\begin{equation}
F_{Y_{n_{1}},\ldots Y_{n_{k}}}\left(  y_{1},\ldots,y_{k}\right)  =\Pr\left\{
\left(  Y_{n_{1}}\leq y_{1}\right)  \wedge\left(  Y_{n_{2}}\leq y_{2}\right)
\wedge\cdots\wedge\left(  Y_{n_{k}}\leq y_{k}\right)  \right\}  .
\end{equation}

For two sequences $a_{m}$ and $b_{m}$, let $a_{m}\sim b_{m}$ denote
$\lim_{m\rightarrow\infty}a_{m}/b_{m}=1$. Let $\operatorname*{const}$ be a
generic positive constant independent of $m$; that is, $\operatorname*{const}$
can have a different value every time it is used. Now $a_{m}=O\left(
b_{m}\right)  $ can be written as $|a_{m}|\leq\operatorname*{const}b_{m}$.

\section{JOINT CUMULATIVE DISTRIBUTION FUNCTION OF ORDER STATISTICS}

First consider the distribution of the order statistics of a random sample
where each sample member is taken from a possibly different population with
its own distribution.

\begin{theorem}
[\citet{Bapat-1989-OSN}, Theorem 4.2]\label{thm:BB}The cumulative distribution
function of the order statistics satisfies
\begin{equation}
F_{Y_{n_{1}},\ldots Y_{n_{k}}}\left(  y_{1},\ldots,y_{k}\right)
=\sum_{\mathbf{i\in}\mathcal{I}}\frac{P_{i_{1},\ldots,i_{k}}\left(
y_{1},\ldots,y_{k}\right)  }{\left(  i_{1}-i_{0}\right)  !\left(  i_{2}%
-i_{1}\right)  !\cdots\left(  i_{k+1}-i_{k}\right)  !}, \label{eq:first}%
\end{equation}
where
\begin{align}
&  P_{i_{1},\ldots,i_{k}}\left(  y_{1},\ldots,y_{k}\right) \nonumber\\
&  \quad=\operatorname{per}%
\begin{bmatrix}
\left[  F_{i}(y_{j})-F_{i}(y_{j-1})\right]  _{\left(  i_{j}-i_{j-1}\right)
\times1}%
\end{bmatrix}
_{j=1,i=1}^{j=k,i=m} \label{eq:BP-perm}%
\end{align}
is the permanent of the block matrix with the block row index $j$ and block
column index $i$. The blocks have $\left(  i_{j}-i_{j-1}\right)  $ rows, and
$1$ column each, which is denoted by the subscript $\left(  i_{j}%
-i_{j-1}\right)  \times1$. Each block has only one distinct entry, which is
$\left[  F_{i}(y_{j})-F_{i}(y_{j-1})\right]  $. We take $F_{i}\left(
y_{0}\right)  =0,\quad F_{i}\left(  y_{k+1}\right)  =1$.

$\qquad$

In expanded form, the permanent (\ref{eq:BP-perm})\ can be written as%
\singlespacing
\begin{equation}
\operatorname{per}%
\begin{bmatrix}
F_{1}(y_{1}) & F_{2}(y_{1}) & \cdots & F_{m}(y_{1})\\
\vdots & \vdots &  & \vdots\\
F_{1}(y_{1}) & F_{2}(y_{1}) & \cdots & F_{m}(y_{1})\\
-----\, & ---- & - & ---\\
F_{1}(y_{2})-F_{1}(y_{1}) & F_{2}(y_{2})-F_{2}(y_{1}) & \cdots & F_{m}%
(y_{2})-F_{m}(y_{1})\\
\vdots & \vdots &  & \vdots\\
F_{1}(y_{2})-F_{1}(y_{1}) & F_{2}(y_{2})-F_{2}(y_{1}) & \cdots & F_{m}%
(y_{2})-F_{m}(y_{1})\\
---- & ---- & - & ----\\
\vdots & \vdots &  & \vdots\\
---- & ---- & - & ----\\
F_{1}(y_{k})-F_{1}(y_{k-1}) & F_{2}(y_{k})-F_{2}(y_{k-1}) & \cdots &
F_{m}(y_{k})-F_{m}(y_{k-1})\\
\vdots & \vdots &  & \vdots\\
F_{1}(y_{k})-F_{1}(y_{k-1}) & F_{2}(y_{k})-F_{2}(y_{k-1}) &  & F_{m}%
(y_{k})-F_{m}(y_{k-1})\\
---- & ---- & - & ----\\
\left[  1-F_{1}\left(  y_{k}\right)  \right]  & \left[  1-F_{2}\left(
y_{k}\right)  \right]  & \cdots & \left[  1-F_{m}\left(  y_{k}\right)  \right]
\\
\vdots & \vdots &  & \vdots\\
\left[  1-F_{1}\left(  y_{k}\right)  \right]  & \left[  1-F_{2}\left(
y_{k}\right)  \right]  & \cdots & \left[  1-F_{m}\left(  y_{k}\right)
\right]
\end{bmatrix}
,
\end{equation}
\jmdouble
where the $j$-th group, $j=1,\ldots,k+1$, contains $i_{j}-i_{j-1}$ repetitions
of the same row.
\end{theorem}

\begin{proof}
The theorem is stated, but not proved in \citet{Bapat-1989-OSN}. We provide a
proof for the sake of completeness, and to prepare the ground for our result.

Define $y_{0}=-\infty$, and $y_{k+1}=\infty$. Note that for $i\in\left\{
1,2,\ldots,m\right\}  $, $F_{i}\left(  y_{0}\right)  =0$, and $F_{i}\left(
y_{k+1}\right)  =1$, since the $F_{i}$ are cumulative distribution functions.
Denote $A=F_{Y_{n_{1}},\ldots Y_{n_{k}}}\left(  y_{1},\ldots,y_{k}\right)
$.\ Then we have%
\begin{equation}
A=\Pr\left(  \bigcap_{j=1}^{k}\left\{  Y_{n_{j}}\leq y_{j}\right\}  \right)
=\Pr\left(  \bigcap_{j=1}^{k}\left\{  \text{at least }n_{j}\text{ of }%
X_{i}\leq y_{j}\right\}  \right) .
\end{equation}
Denote by $I_{j}$ the random variable equal to the number of $X_{i}$ such that
$X_{i}\leq y_{j}$. Then $I_{1}\leq I_{2}\leq\cdots\leq I_{k}$, and the
condition that at least $n_{j}$ of $X_{i}\leq y_{j}$ is equivalent to
$I_{j}\geq n_{j}$. Thus,
\begin{equation}
A=\Pr\left(  \bigcap_{j=1}^{k}\left\{  I_{j}\geq n_{j}\right\}  \right)
=\Pr\left(  \bigcup_{\mathbf{i\in}\mathcal{I}}^{i_{2}}\bigcap_{j=1}%
^{k}\left\{  I_{j}=i_{j}\right\}  \right)  ,
\end{equation}
and, since the events $\bigcap_{j=1}^{k}\left\{  I_{j}=i_{j}\right\}  $ for
different $\mathbf{i}$ are disjoint,
\begin{align}
A  &  =\sum_{\mathbf{i\in}\mathcal{I}}\Pr\left(  \bigcap_{j=1}^{k}\text{
}\left\{  I_{j}=i_{j}\right\}  \right) \label{eq:penultimate}\\
&  =\sum_{\mathbf{i\in}\mathcal{I}}\Pr\left(  \bigcap_{j=1}^{k+1}\text{
}\left\{  \text{exactly }i_{j}-i_{j-1}\text{ of }X_{i}\in(y_{j-1}%
,y_{j}]\right\}  \right)  . \label{eq:ultimate}%
\end{align}
\qquad

Now fix $\mathbf{i}$ and write an arbitrary permutation of $\left\{
1,2,\ldots,m\right\}  $ as
\begin{equation}
\pi=\left(  \pi_{1},\pi_{2},\ldots,\pi_{k},\pi_{k+1}\right)  ,
\end{equation}
where each subsequence $\pi_{j}$ has exactly $i_{j}-i_{j-1}$ terms.  We will
use $\left\{  \pi_{j}\right\}  $ to denote the set of the terms. Then,
\begin{align}
\text{ }  &  \exists\pi\forall j\in\left\{  1,2,\ldots,k+1\right\}
:\text{exactly }i_{j}-i_{j-1}\text{ of }X_{i}\in(y_{j-1},y_{j}]\quad\\
&  \Longleftrightarrow\exists\pi\forall j\in\left\{  1,2,\ldots,k+1\right\}  :
\forall i\in\left\{  \pi_{j}\right\}  :X_{i}\in(y_{j-1},y_{j}].
\end{align}
Hence,%
\begin{align}
&  \Pr\left(  \bigcap_{j=1}^{k+1}\text{ }\left\{  \text{exactly }i_{j}%
-i_{j-1}\text{ of }X_{i}\in(y_{j-1},y_{j}]\right\}  \right) \\
&  =\frac{\sum_{\pi}\Pr\left(  \bigcap_{j=1}^{k+1}\bigcap_{i\in\left\{
\pi_{j}\right\}  }\text{ }\left\{  X_{i}\in(y_{j-1},y_{j}]\right\}  \right)
}{\left(  i_{1}-i_{0}\right)  !\cdots\left(  i_{k+1}-i_{k}\right)  !}\\
&  =\frac{\sum_{\pi}\prod\limits_{j=1}^{k+1}\prod\limits_{i\in\left\{  \pi
_{j}\right\}  }\left[  F_{i}\left(  y_{j}\right)  -F_{i}\left(  y_{j-1}%
\right)  \right]  }{\left(  i_{1}-i_{0}\right)  !\cdots\left(  i_{k+1}%
-i_{k}\right)  !},
\end{align}
because the events in the intersection are independent:\ there is one event
for each $X_{i}$, which are independent random variables. Substituting into
(\ref{eq:penultimate})\ and comparing with the definition of the permanent
(\ref{eq:permanent}) concludes the proof.
\end{proof}

As noted in the introduction, using a general algorithm for permanents is
prohibitively expensive. Given simplifying assumptions, however, the problem
becomes easier. In the case when the variables $X_{1}$, $X_{2},\ldots,X_{m}$
are independent and identically distributed (that is, the classical case of
sampling from a single population), Theorem \ref{thm:BB} reduces to the
following well-known result {\citep[p. 11]{David-2003-OS}}.

\begin{theorem}
\label{thm:one-population}Suppose that $F_{i}=F$ for all $i$. Then the joint
cumulative distribution function of the order statistics satisfies
\begin{equation}
F_{Y_{n_{1}},\ldots Y_{n_{k}}}\left(  y_{1},\ldots,y_{k}\right)
=\sum_{\mathbf{i\in}\mathcal{I}}m!\prod\limits_{j=1}^{k+1}\frac{\left[
F\left(  y_{j}\right)  -F\left(  y_{j-1}\right)  \right]  ^{i_{j}-i_{j-1}}%
}{\left(  i_{j}-i_{j-1}\right)  !}.
\end{equation}

\end{theorem}

Now consider drawing a random sample from two populations, each with a
different cumulative distribution function, say $F\left(  x\right)  $, and
$G\left(  x\right)  $. Sample the first $n$ random variables from the first
population with the distribution function $F$, and then $m-n$ from the second
population with the distribution function $G$. Then the permanents from
Equation \ref{eq:first} (\citet{Bapat-1989-OSN}) simplify to the block form
with constant blocks,%
\begin{align}
&  P_{i_{1},\ldots,i_{k}}\left(  y_{1},\ldots,y_{k}\right) \nonumber\\
&  =\operatorname{per}%
\begin{bmatrix}
\left[  F(y_{1})-F(y_{0})\right]  _{\left(  i_{1}-i_{0}\right)  \times n} &
\left[  G(y_{1})-G(y_{0})\right]  _{\left(  i_{1}-i_{0}\right)  \times\left(
m-n\right)  }\\
\left[  F(y_{2})-F(y_{1})\right]  _{\left(  i_{2}-i_{1}\right)  \times n} &
\left[  G(y_{2})-G(y_{1})\right]  _{\left(  i_{2}-i_{1}\right)  \times\left(
m-n\right)  }\\
\vdots & \vdots\\
\left[  F(y_{k+1})-F(y_{k})\right]  _{\left(  i_{k+1}-i_{k}\right)  \times n}
& \left[  G(y_{k+1})-G(y_{k})\right]  _{\left(  i_{k+1}-i_{k}\right)
\times\left(  m-n\right)  }%
\end{bmatrix}
, \label{eq:2-perm}%
\end{align}
where the subscripts indicate the dimensions of blocks created by the
repetition of the term in the brackets, and we take
\begin{equation}
F\left(  y_{0}\right)  =G\left(  y_{0}\right)  =0,\quad F\left(
y_{k+1}\right)  =G\left(  y_{k+1}\right)  =1. \label{eq:FG-added}%
\end{equation}

In expanded form, the permanent (\ref{eq:2-perm})\ can be written as%
\singlespacing
\begin{equation}
\text{per}%
\begin{bmatrix}
F(y_{1}) & \cdots & F(y_{1}) & G(y_{1}) & \cdots & G(y_{1})\\
\vdots &  & \vdots & \vdots &  & \vdots\\
F(y_{1}) & \cdots & F(y_{1}) & G(y_{1}) & \cdots & G(y_{1})\\
-----\, & - & ---- & --- & - & ---\\
F(y_{2})-F(y_{1}) & \cdots & F(y_{2})-F(y_{1}) & G(y_{2})-G(y_{1}) & \cdots &
G(y_{2})-G(y_{1})\\
\vdots &  & \vdots & \vdots &  & \\
F(y_{2})-F(y_{1}) & \cdots & F(y_{2})-F(y_{1}) & G(y_{2})-G(y_{1}) & \cdots &
G(y_{2})-G(y_{1})\\
---- & - & ---- & ---- & -\, & ---\\
\vdots &  & \vdots & \vdots &  & \vdots\\
---- & - & ---- & ---- & - & ---\\
F(y_{k})-F(y_{k-1}) & \cdots & F(y_{k})-F(y_{k-1}) & G(y_{k})-G(y_{k-1}) &
\cdots & G(y_{k})-G(y_{k-1})\\
\vdots &  & \vdots & \vdots &  & \vdots\\
F(y_{k})-F(y_{k-1}) & \cdots & F(y_{k})-F(y_{k-1}) & G(y_{k})-G(y_{k-1}) &
\cdots & G(y_{k})-G(y_{k-1})\\
---- & - & ---- & ---- & - & ---\\
1-F\left(  y_{k}\right)  & \cdots & 1-F\left(  y_{k}\right)  & 1-G\left(
y_{k}\right)  & \cdots & 1-G\left(  y_{k}\right) \\
\vdots &  & \vdots & \vdots &  & \vdots\\
1-F\left(  y_{k}\right)  & \cdots & 1-F\left(  y_{k}\right)  & 1-G\left(
y_{k}\right)  & \cdots & 1-G\left(  y_{k}\right)
\end{bmatrix}
.
\end{equation}
\jmdouble

This special form of the permanent allows us to evaluate the joint
distribution of the order statistic more efficiently.

\begin{theorem}
\label{thm:two-polulations}Suppose that $F_{i}\left(  x\right)  =F\left(
x\right)  $, for all $1\leq i\leq n,$ and $F_{i}\left(  x\right)  =G\left(
x\right)  $, for all $n+1\leq i\leq m$. Then
\begin{align}
&  F_{Y_{n_{1}},\ldots Y_{n_{k}}}\left(  y_{1},\ldots,y_{k}\right)
=\nonumber\\
&  \sum_{\mathbf{i\in}\mathcal{I}}\sum_{\boldsymbol{\lambda}}\prod_{j=1}%
^{k+1}\frac{n!\left(  m-n\right)  !}{\lambda_{j}!\left(  i_{j}-i_{j-1}%
-\lambda_{j}\right)  !}\nonumber\\
&  \qquad\cdot\left[  F\left(  y_{j}\right)  -F\left(  y_{j-1}\right)
\right]  ^{\lambda_{j}}\left[  G\left(  y_{j}\right)  -G\left(  y_{j-1}%
\right)  \right]  ^{i_{j}-i_{j-1}-\lambda_{j}},
\end{align}
where $\boldsymbol{\lambda}=\left(  \lambda_{1},\lambda_{2},\ldots
,\lambda_{k+1}\right)  $ ranges over all integer vectors such that%
\begin{equation}
\lambda_{1}+\lambda_{2}+\cdots+\lambda_{k+1}=n,\quad0\leq\lambda_{j}\leq
i_{j}-i_{j-1}. \label{eq:lambda-constr}%
\end{equation}

\end{theorem}

\begin{proof}
We evaluate the permanents $P_{i_{1},\ldots,i_{k}}\left(  y_{1},\ldots
,y_{k}\right)  $ from (\ref{eq:2-perm}). Let $S_{1}=\left\{  1,2,\ldots
,n\right\}  $ and $S_{2}=\left\{  n+1,n+2,\ldots,m\right\}  $. Write a
permutation of $\left\{  1,2,\ldots,m\right\}  $ as $\pi=\left(  \pi_{1}%
,\pi_{2},\ldots,\pi_{k},\pi_{k+1}\right)  $, where each subsequence $\pi_{j}$
has exactly $i_{j}-i_{j-1}$ terms. The subsequence $\pi_{j}$ is a list of the subscripts of the random variables that fall in the interval $\left(y_{j-1},y_{j}\right)$. Then the term in the definition of the
permanent (\ref{eq:permanent}) associated with $\pi$ is%
\begin{equation}
\prod_{i=1}^{m}a_{i,\pi\left(  i\right)  }=\prod_{j=1}^{k+1}\left[  F\left(
y_{j}\right)  -F\left(  y_{j-1}\right)  \right]  ^{\lambda_{j}}\left[
G\left(  y_{j}\right)  -G\left(  y_{j-1}\right)  \right]  ^{i_{j}%
-i_{j-1}-\lambda_{j}},
\end{equation}
where $\lambda_{j}$ is the number of random variables with subscripts listed in  $\left\{  \pi
_{j}\right\}  $ that are in $S_{1}$. For illustration, the intervals and the number of order statistics of each
type in them are shown in Table \ref{tab:2-pop}.

\begin{table}[ptb]%
\begin{tabular}
[c]{c|cccc|c}%
Interval & $(-\infty,y_{1}]$ & $(y_{1},y_{2}]$ & $\cdots$ & $\left(
y_{k},\infty\right)  $ & Total\\\hline
$\text{\# $\in S_{1}$\ }$ & $\lambda_{1}$ & $\lambda_{2}$ & $\cdots$ &
$\lambda_{k+1}$ & $n$\\
$\text{\# $\in S_{2}$ \ }$ & $i_{1}-\lambda_{1}$ & $i_{2}-i_{1}-\lambda_{2}$ &
$\cdots$ & $m-i_{k}-\lambda_{k+1}$ & $m-n$\\\hline
Total & $i_{1}$ & $i_{2}-i_{1}$ & $\cdots$ & $m-i_{k}$ & $m$%
\end{tabular}
\caption{Total number of order statistics in each interval, and number from
population 1 and 2 in each interval.}%
\label{tab:2-pop}%
\end{table}

The number of permutations $\pi$ such
that $\lambda_{j}$ is the number of the elements from $\left\{  \pi
_{j}\right\}  $ that are in $S_{1}$ is found as the product $ABC$, where%
\begin{equation}
A=\frac{n!}{\prod_{j=1}^{k+1}\lambda_{j}!}%
\end{equation}
is the number of ways to distribute the $n$ elements of $S_{1}$ so that set
$j$ has $\lambda_{j}$ elements (the multinomial coefficient),
\begin{equation}
B=\frac{\left(  m-n\right)  !}{\prod_{j=1}^{k+1}\left(  i_{j}-i_{j-1}%
-\lambda_{j}\right)  !}%
\end{equation}
is the number of ways to distribute the $m-n$ elements of $S_{1}$ so that set
$j$ has $i_{j}-i_{j-1}-\lambda_{j}$ elements, and%
\begin{equation}
C=\prod_{j=1}^{k+1}\left(  i_{j}-i_{j-1}\right)  !
\end{equation}
is the number of permutations that do not change the distribution of the
elements $S_{1}$ and $S_{2}$ into those sets. Thus,%
\begin{align}
&  P_{i_{1},\ldots,i_{k}}\left(  y_{1},\ldots,y_{k}\right)  =\sum_{\pi}%
\prod_{i=1}^{m}a_{i,\pi\left(  i\right)  }\nonumber\\
&  =\sum_{\boldsymbol{\lambda}}\prod_{j=1}^{k+1}\frac{\left(  i_{j}%
-i_{j-1}\right)  !}{\lambda_{j}!\left(  i_{j}-i_{j-1}-\lambda_{j}\right)
!}\nonumber\\
&  \qquad\cdot\left[  F\left(  y_{j}\right)  -F\left(  y_{j-1}\right)
\right]  ^{\lambda_{j}}\left[  G\left(  y_{j}\right)  -G\left(  y_{j-1}%
\right)  \right]  ^{i_{j}-i_{j-1}-\lambda_{j}},
\end{align}
with the sum over all $\boldsymbol{\lambda}$ that satisfy
(\ref{eq:lambda-constr}). The result now follows from Theorem~\ref{thm:BB}.
\end{proof}

The proof of Theorem \ref{thm:two-polulations} easily carries over to the
general case of order statistics of a sample selected from an arbitrary number
of populations. The proof of the next theorem can therefore be omitted.

\begin{theorem}
\label{thm:arbitrary}Suppose that $F_{i}=G_{1}$ for the first $m_{1}$ indices
$i$, $F_{i}=G_{2}$ for the next $m_{2}$ indices $i$, etc., and $F_{i}=G_{N}$
for the last $m_{N}$ indices $i$, with%
\begin{equation}
m_{1}+\cdots+m_{N}=m,\quad m_{s}>0\text{ for all }s.
\end{equation}
Then
\begin{align}
&  F_{Y_{n_{1}},\ldots Y_{n_{k}}}\left(  y_{1},\ldots,y_{k}\right)  =\\
&  =\sum_{\mathbf{i\in}\mathcal{I}}\sum_{\left[  \lambda_{js}\right]  }%
\prod_{j=1}^{k+1}\prod\limits_{s=1}^{N}\frac{m_{s}!}{\lambda_{js}!}\left[
G_{s}\left(  y_{j}\right)  -G_{s}\left(  y_{j-1}\right)  \right]
^{\lambda_{js}}%
\end{align}
where the summation is over all integer matrices $\left[  \lambda_{js}\right]
$ size $k+1$ by $N$ such that%
\begin{align}
\lambda_{js}  &  \geq0\quad\text{for all }j\text{ and all }s,\\
\sum_{j=1}^{k+1}\lambda_{js}  &  =m\quad\text{for all }s,\\
\sum_{s=1}^{N}\lambda_{js}  &  =i_{j}-i_{j-1}\quad\text{for all }j,
\end{align}
and we take $G_{s}\left(  y_{0}\right)  =0$, $G_{s}\left(  y_{k+1}\right)  =1$.
\end{theorem}

Theorem \ref{thm:arbitrary} covers all of the theorems above. In the
particular case when all $m_{i}=1$, i.e., every distribution is different
because it comes from a different population, it gives exactly the same result
as Theorem \ref{thm:BB}. With two populations, the complexity of
Theorem \ref{thm:arbitrary} reduces to the complexity of Theorem
\ref{thm:two-polulations}. The complexity of Theorem \ref{thm:two-polulations}
is less than that of the Theorem \ref{thm:BB} from \citet{Bapat-1989-OSN}, as
discussed in the next section.

\section{COMPLEXITY}

We will now compare the relative complexity of Theorem \ref{thm:BB}, from
\citet{Bapat-1989-OSN}, and our formula, Theorem \ref{thm:two-polulations}. We
assume that the evaluation of the cumulative distribution function of each of
the statistics takes a constant number of operations.

For $1\leq n_{1}<n_{2}<\cdots<n_{k}\leq m$, denote the number of elements of
the index set $\mathcal{I}$ by
\begin{equation}
\nu\left(  n_{1},n_{2},\cdots,n_{k};m\right)  =\left\vert \mathcal{I}%
\right\vert =\sum_{i_{k}=n_{k}}^{m}\sum_{i_{k-1}=n_{k-1}}^{i_{k}}\cdots
\sum_{i_{1}=n_{1}}^{i_{2}}1. \label{eq:summation}%
\end{equation}

\begin{theorem}
\label{lem:catalan-triangle}The number $\nu\left(  n_{1},n_{2},\cdots
,n_{k};m\right)  $ of the Bapat-Beg permanents in Theorem \ref{thm:BB} is
bounded by%
\begin{equation}
\nu\left(  n_{1},n_{2},\ldots,n_{k};m\right)  \leq\nu\left(  1,2,\ldots
,k;m\right)  \leq\nu\left(  1,2,\ldots,m;m\right)  =C_{m}, \label{eq:I-ineq}%
\end{equation}
where%
\begin{equation}
\nu\left(  1,2,\ldots,k;m\right)  =\binom{m+k}{k}\left(  1-\frac{k}%
{m+1}\right)  , \label{eq:catalan-triangle}%
\end{equation}
and
\begin{equation}
C_{m}=\frac{1}{m+1}\binom{2m}{m}=\frac{(2m)!}{(m+1)!\,m!}. \label{eq:catalan}%
\end{equation}

\end{theorem}

\begin{proof}
The inequalities in (\ref{eq:I-ineq}) are obtained by taking the smallest
numbers for $n_{1},n_{2},\ldots,n_{k}$ and the largest possible value for $k$,
which both give the largest number of terms. We now prove that%
\begin{equation}
\nu\left(  1,2,\ldots,k;m\right)  =\binom{m+k}{k}-\binom{m+k}{k-1}
\label{eq:induction}%
\end{equation}
by induction over $k$. For $k=1$, (\ref{eq:induction}) follows from%
\begin{equation}
\nu\left(  1;m\right)  =\sum_{i_{1}=1}^{m}1=m
\end{equation}
and%
\begin{equation}
\binom{m+1}{1}-\binom{m+1}{1-1}=\left(  m+1\right)  -1=m.
\end{equation}

Now assume that (\ref{eq:induction}) holds for some $k$ and we will show that
\begin{equation}
\nu\left(  1,2,\ldots,k+1,m\right)  =\binom{m+k+1}{k+1}-\binom{m+k+1}{k}.
\label{eq:k-plus-1}%
\end{equation}
From the definition (\ref{eq:summation}) and the induction assumption
(\ref{eq:induction}), it follows that%
\begin{align}
\nu\left(  1,2,\ldots,k+1;m\right)   &  =\sum_{i_{k+1}=k+1}^{m}\nu\left(
1,2,\ldots,k;i_{k+1}\right) \\
&  =\sum_{i=k+1}^{m}\binom{i+k}{k}-\binom{i+k}{k-1}\\
&  =\sum_{i=k+1}^{m}\left[  \binom{i+k+1}{k+1}-\binom{i+k}{k+1}\right] \\
&  \quad-\sum_{i=k+1}^{m}\left[  \binom{i+k}{k}-\binom{i+k+1}{k}\right]  ,
\end{align}
where we have used the identity
\begin{equation}
\binom{n}{r}-\binom{n-1}{r}=\binom{n-1}{r-1}%
\end{equation}
twice. Both sums telescope, and we get
\begin{align}
\nu\left(  1,2,\ldots,k+1;m\right)   &  =\left[  \binom{m+k+1}{k+1}%
-\binom{2k+1}{k+1}\right] \\
&  -\left[  \binom{m+k+1}{k}+\binom{2k+1}{k}\right]  ,
\end{align}
which, noting that%
\begin{equation}
\binom{2k+1}{k+1}=\frac{\left(  2k+1\right)  !}{\left(  k+1\right)
!k!}=\binom{2k+1}{k},
\end{equation}
gives (\ref{eq:k-plus-1}). Equations (\ref{eq:catalan-triangle})\ and
(\ref{eq:catalan}) follow from (\ref{eq:induction}) by a direct computation:
\begin{align}
\binom{m+k}{k}-\binom{m+k}{k-1}  &  =\frac{m+k}{1}\frac{m+k-1}{2}\cdots
\frac{m+2}{k-1}\frac{m+1}{k}\\
&  -\frac{m+k}{1}\frac{m+k-1}{2}\cdots\frac{m+2}{k-1}\\
&  =\binom{m+k}{k}\left(  1-\frac{k}{m+1}\right)  ,
\end{align}
and%
\begin{equation}
\binom{m+m}{m}-\binom{m+m}{m-1}=\binom{2m}{m}\left(  1-\frac{m}{m+1}\right)
=\frac{1}{m+1}\binom{2m}{m},
\end{equation}
which concludes the proof.
\end{proof}

The numbers $C_{m}$ defined by (\ref{eq:catalan}) are known as the Catalan
numbers \citep{Stanley-1999-EC}, and the numbers $a_{k,m}=\nu\left(
1,2,\ldots,k;m\right)  $ are called the Catalan triangle
\citep{Shapiro-1976-CT}. From the Stirling approximation $m!\sim\sqrt{2\pi
m}\,m^{m}/e^{m}$, the growth of Catalan numbers is exponential,%

\begin{equation}
C_{m}\sim\operatorname*{const}m^{-3/2}4^{m}>\operatorname*{const}\alpha
^{m},\quad\label{eq:catalan-exponential}%
\end{equation}
for any $1<\alpha<4$ (with a different $\operatorname*{const}$ for each
$\alpha$).

\begin{theorem}
The worst case complexity of computing the distribution function of the order
statistics from Theorem \ref{thm:BB} is%
\begin{equation}
\operatorname*{const}C_{m}mK^{m}\sim\operatorname*{const}m^{-1/2}4^{m} P(m),
\label{eq:complexity-BB}%
\end{equation}
where $P(m)$ is the number of operations for computing permanent of order $m$.
\end{theorem}

\begin{proof}
The denominator in (\ref{eq:first}) requires at most $O\left(  m\right)  $
operations, and there are at most $C_{m}$ terms in the sum by Theorem
\ref{lem:catalan-triangle}.
\end{proof}

It is known that the complexity of computing the permanent is bounded by
\[
P(m)=O(m^{a}2^{m})
\]
for some $a$, e.g., from the Ryser's formula \citep{Knuth-1998-TAC}. So, the
complexity of the computation of the distribution function from Theorem
\ref{thm:BB} is exponential in $m$. Therefore, the computation is practical
only for small $m$.

Fortunately, a drastic reduction of complexity is possible in the case when
the order statistics come from two populations. In fact, the
complexity reduces still farther when we need only a small number $k$ of order statistics.

\begin{theorem}
\label{thm:two-k}Let $C\left(  k,m,n\right)  $ be the number of operations  in
Theorem \ref{thm:two-polulations} to evaluate the joint distribution function
of $k$ order statistics from $m$ random variables from two populations, with
$n\leq m$ of the variables from the first population. Then%
\begin{equation}
C\left(  k,m,n\right)  \leq\operatorname*{const}k\binom{m+k}{k}\binom{n+k}%
{k}\left(  1-\frac{k}{m+1}\right)  .\label{eq:two-k}%
\end{equation}
In the worst case over all $k$ and $n$, the complexity is bounded by
\begin{equation}
C\left(  k,m,n\right)  \leq\operatorname*{const}m\frac{\left(  2m\right)
^{2}}{\left(  m!\right)  ^{4}}\sim\operatorname*{const}16^{m}%
,\label{eq:two-k-worst-case}%
\end{equation}
For any fixed $k$, the complexity is bounded by%
\begin{equation}
~C\left(  k,m,n\right)  =O\left(  m^{k}n^{k}\right)  .\label{eq:two-fixed-k}%
\end{equation}
i.e., the complexity is polynomial in $m$.
\end{theorem}

\begin{proof}
The complexity is bounded by $\operatorname*{const}CLM$, where $C=\binom
{m+k}{k}\left(  1-\frac{k}{m}\right)  $ is the number of terms in the sum over
$\mathbf{i}$, $L$ is the number of possible index vectors
$\boldsymbol{\lambda}$ satisfying (\ref{eq:lambda-constr}), and $M$ is the
complexity of evaluating the products in one term of the sum, which is
$M=O\left(  k\right)  $. To bound $L$, drop the upper bounds in
(\ref{eq:lambda-constr}). Thus $L$ is bounded above by the number of all
integer vectors $\boldsymbol{\lambda}$ such that%
\begin{equation}
\lambda_{1}+\lambda_{2}+\cdots+\lambda_{k+1}=n,\quad\lambda_{j}\geq0\text{ for
all }j,
\end{equation}
which is the same as the number of ways to distribute $n$ indistinguishable
objects to $k+1$ distinguishable bins, which equals to $\binom{n+k}{k}$. This
gives (\ref{eq:two-k}).

The bound (\ref{eq:two-k-worst-case}) follows by taking a pessimistic value of
$k$ in each term (\ref{eq:two-k}) - twice $k=m$, then $k=0$, and pessimistic
value $n=m$. The second part of (\ref{eq:two-k-worst-case}) follows from the
Stirling formula.

The polynomial bound (\ref{eq:two-fixed-k}) follows from (\ref{eq:two-k}) and
the inequality%
\[
\binom{p+k}{k}=\frac{\left(  p+k\right)  \left(  p+k-1\right)  \cdots\left(
p+1\right)  }{1\cdot2\cdots k}\leq\operatorname*{const}(k)p^{k}%
\]
applied with $p=m$ and $p=n$.
\end{proof}

Although the complexity of evaluating the cumulative distribution function of
order statistics from Theorem \ref{thm:BB} is exponential in the general case,
we have shown in Theorem \ref{thm:two-k} that the complexity is bounded by a
polynomial of a small degree when there are only two populations, and the number
of order statistics considered, $k$, is fixed and small. The complexity also depends on $n$, the number of random variables from the first population, $S_{1}$.  In general, $n$ is fixed by the state of nature.

\begin{figure}[ptb]
\begin{center}
\includegraphics[width=3.5in]{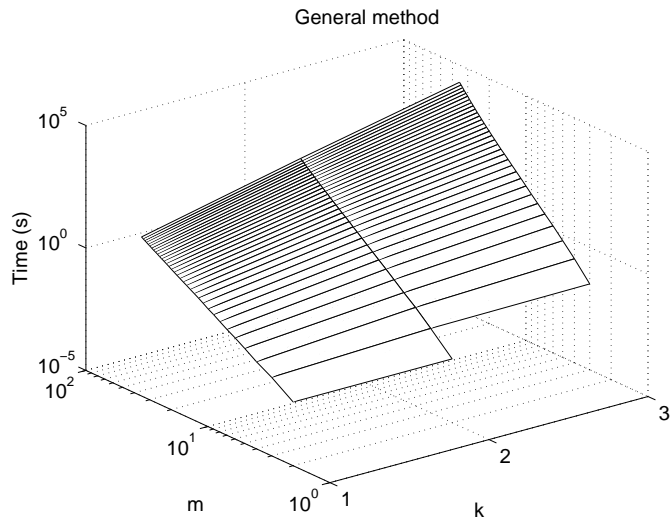}
\end{center}
\caption{Times for evaluating the joint cumulative distribution function of
the first $k$ order statistics of $m$ random variables from two distributions,
using the general Bapat-Beg formula (Theorem \ref{thm:BB}).}%
\label{fig:bb}%
\end{figure}

\begin{figure}[ptb]
\begin{center}
\includegraphics[width=3.5in]{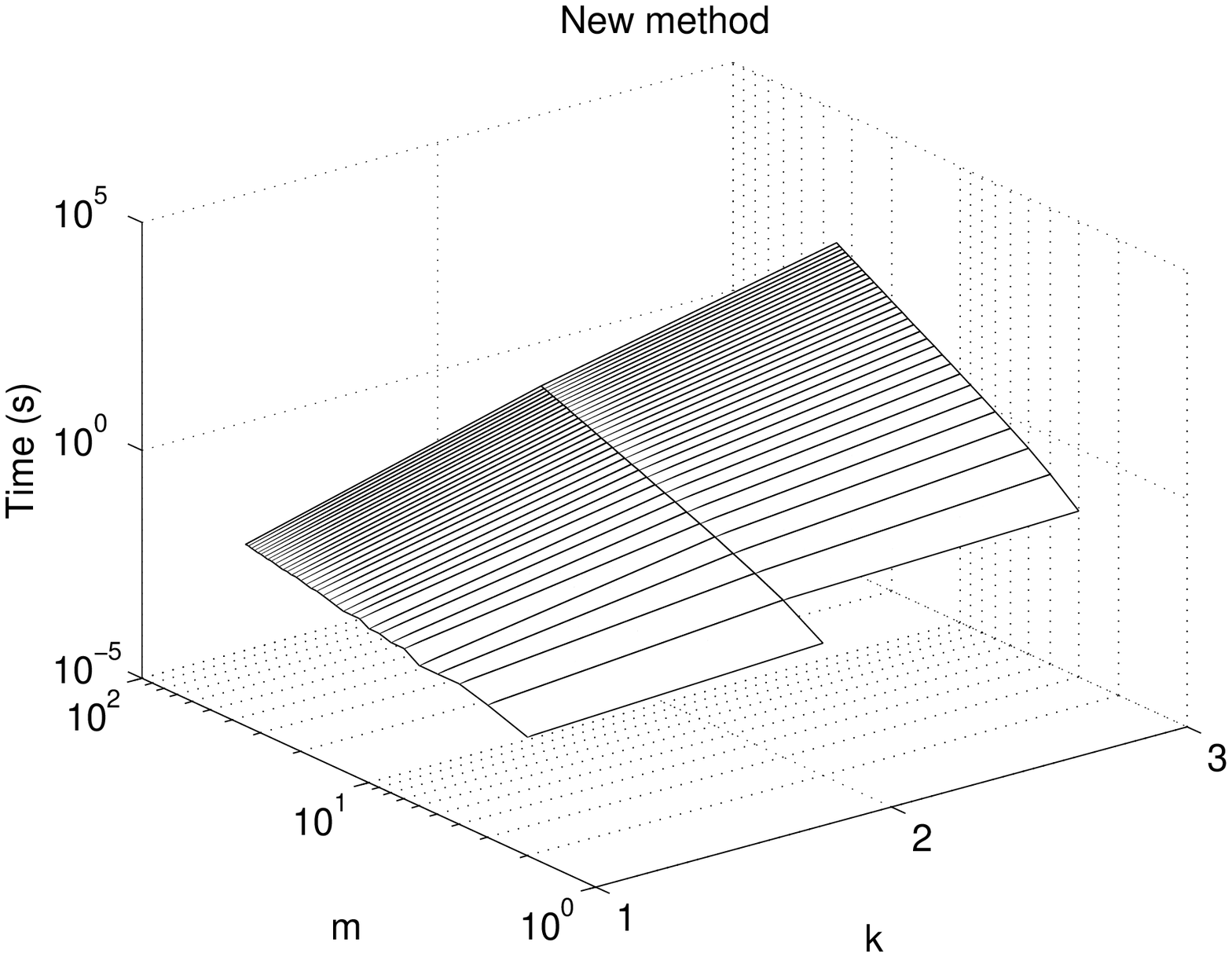}
\end{center}
\caption{Times for evaluating the joint cumulative distribution function of
the first $k$ order statistics of $m$ random variables from two distributions,
using the new formula from Theorem \ref{thm:two-polulations}.}%
\label{fig:nw}%
\end{figure}

\begin{figure}[ptb]
\begin{center}
\includegraphics[width=3.5in]{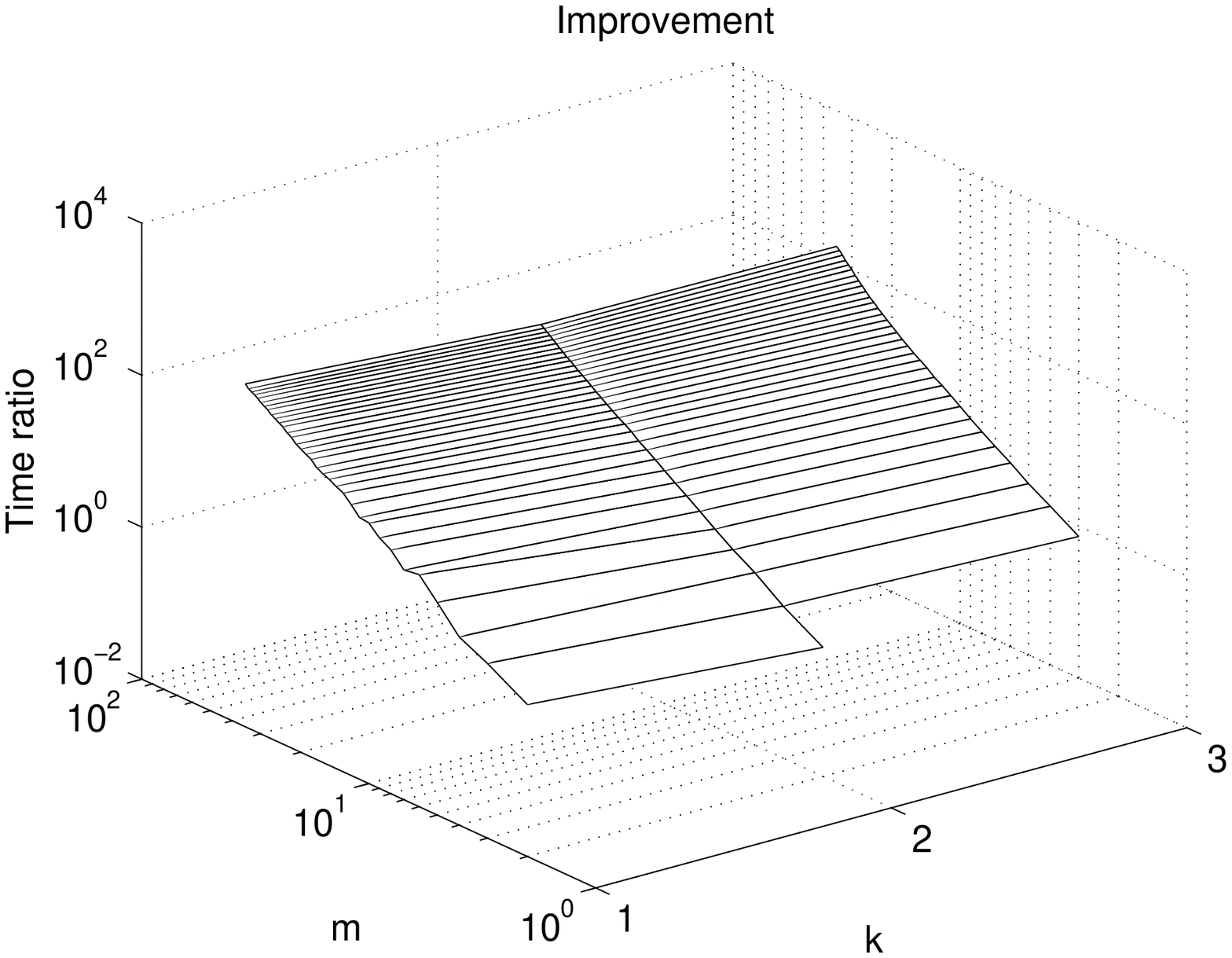}
\end{center}
\caption{Ratio of times for evaluating the joint cumulative distribution
function of the first $k$ order statistics of $m$ random variables from two
distributions, using the Bapat-Beg formula (Theorem \ref{thm:BB}) and the new
formula from Theorem \ref{thm:two-polulations}.}%
\label{fig:im}%
\end{figure}

\begin{table}[ptb]%
\begin{tabular}
[c]{l|l|l}%
Bapat-Beg formula & New formula & Improvement\\
Theorem \ref{thm:BB}, Fig.~\ref{fig:bb} & Theorem \ref{thm:two-polulations},
Fig.~\ref{fig:nw} & Fig.~\ref{fig:im}\\\hline
&  & \\
$10^{-2.9-0.36k}m^{2.0+1.1k}$ & $10^{-2.6-0.01k}m^{0.06+1.02k}$ &
$10^{-0.30-0.34k}m^{1.93+0.09k}$%
\end{tabular}
\caption{Fit of timing in Mathematica of the evaluation of the joint
distribution of the first $k$ statistics of $m$ variables from two populations
($n$=1 from one population, $m-n$ from the other). For fixed $k$, regression
was used to fit the logarithm of the time with a linear function of $\log m$,
and regression was then used again to fit the coefficients by linear functions
of $k$.}%
\label{tab:fit}%
\end{table}

To confirm and illustrate the result, we have conducted a timing experiment. We calculated
the joint distribution function in the case of two populations.  We considered $k=1$, $k=2$, and $k=3$, and fixed $n=1$. We measured the amount of time it took to compute the joint distribution function
using the general Bapat Beg formula with permanents (Fig.~\ref{fig:bb}) and
the new special formula (Fig.~\ref{fig:nw}). Both theorems were implemented in
Mathematica . The permanents were computed in Mathematica using the
code
\begin{align*}
\mathtt{Permanent[A\_List]:=} &  \mathtt{With[{v=Array[x,Length[A]]},}\\
&  \mathtt{Coefficient[Times@@(A.v),Times@@v]}%
\end{align*}
from \citet{Weisstein-2007-P}. This function computes the permanent of matrix
$A$ by Vardi's formula as the coefficient of $x_{1}\cdots x_{m}$ in
\[
\prod_{i=1}^{m}\left(  a_{i1}x_{1}+a_{i2}x_{2}+\cdots+a_{im}x_{m}\right)  ,
\]
using symbolic manipulation with automatic caching of partial results by the
Mathematica kernel. Amazingly, calculating the permanent from (\ref{eq:2-perm}%
) in Mathematica results in times that grow polynomially with $m$, the number
of rows in the permanent. Consequently, for two populations, while the
theoretical complexity of Bapat Beg is exponential, the actual time observed
while calculating the formulas in Mathematica was polynomial
(Fig.~\ref{fig:bb}). Graphing the time versus the log of $m$ produces almost
straight lines in a log-log plot. We attribute this speedup to the reuse of
partial results by the Mathematica kernel.

Mathematica calculates the Bapat Beg formula more rapidly than
predicted. In the timing experiment, the observed times for the new formula (Theorem
\ref{thm:two-polulations}) are much faster than the Bapat Beg formula. The observed improvement was quite
dramatic (Fig.~\ref{fig:im}). The observed improvement is of the order $m^{2}$
(Table \ref{tab:fit}). The observed complexity of the new formula for two
populations was of the order $m^{k}$, which confirms the result of Theorem
\ref{thm:two-k} for constant $n=1$.

All calculations were done using a custom New Tech Solutions workstation with
4 AMD Opteron 848 processors running Mathematica 5.2, under the SuSE Linux
Enterprise Server 10 operating system.

Mathematica code to calculate the cumulative distribution function for
arbitrary collections of order statistics of independent random variables
which may have different distributions is available free from the authors.
Examples demonstrating the use of the software are also available from the authors.


\end{document}